	\documentclass[11pt,titlepage]{article}
	\usepackage{amsfonts,emlines}

\newcommand{\qed}[1][6]{\hbox{\vrule width #1pt height #1pt depth 0pt}}
\newtheorem{thm}{Theorem}[section]		
\newtheorem{prop}{Proposition}[section]		
\newtheorem{defi}{Definition}[section]		
\def\openone{\leavevmode\hbox{\small1\kern-3.8pt\normalsize1}}
\title{\vspace*{-1.84in}\bfseries
\vspace*{-7ex}
{
\begin{flushright}
	  \textbf{\large LANL xxx E-print archive No. dg-ga/9709017}\\[5ex]
\end{flushright}
}
\huge Linear transports along paths \\[6pt] in vector bundles \\
\vspace{10pt}
\LARGE V. Properties of curvature and torsion
}
\vspace{10pt}
\renewcommand{\thefootnote}{\fnsymbol{footnote}}
\author{Bozhidar Z. Iliev
\thanks{Permanent address:
Department Mathematical Modeling,
Institute for Nuclear Research and \mbox{Nuclear} Energy,
Bulgarian Academy of Sciences,
Boul. Tzarigradsko chauss\'ee~72, 1784 Sofia, Bulgaria.
}
\thanks{E-mail address: bozho@inrne.bas.bg}
\thanks{URL: http://www.inrne.bas.bg/mathmod/bozhome/}
}

\newlength{\bo}
\newlength{\ho}
\newlength{\up}
\newlength{\down}
\newlength{\middle}
\newcommand{\BOZHO}{\leavevmode\hbox{\slshape\bfseries
\settowidth{\bo}{BO}%
\settowidth{\ho}{HO}%
\settowidth{\middle}{/}%
\settoheight{\up}{BOZHO}%
\settodepth{\down}{/}%
\addtolength{\up}{+0.15\up}%
\addtolength{\bo}{+\middle}%
\rule[\up]{\bo}{0.15ex}%
\hspace{-\bo}BO%
\hspace{+0.09em}\raisebox{+0.17\up}{/}%
\hspace{-0.20em}\raisebox{+0.71\up}{$\bullet$}%
\hspace{-0.33em}\hspace{-1.14\middle}\raisebox{-0.4\up}{$\bullet$}%
\hspace{-0.30em}%
\addtolength{\down}{-0.41\down}%
\addtolength{\ho}{+1.5\middle}%
\rule[-\down]{\ho}{0.15ex}%
\addtolength{\ho}{-\middle}%
\hspace{-\ho}\hspace{+0.18em}%
\raisebox{+0.17\up}{HO}%
}}

\date{
 \vspace{5ex}
 Ended: November 28, 1995  \\[3pt]
 Revised: March 25, November 6, 1996  \\
 Updated: July 9, 1998\\[0.6ex]
 Produced: \today\\
 Submitted to JINR communication: January 13, 1997\\
 Published: Communication JINR, E5-97-1, Dubna, 1997\\[1ex]
 LANL xxx archive E-print No. dg-ga/9709017 \\[5ex]
 {\Huge\BOZHO}
}

\begin{document}

\renewcommand{\thefootnote}{\fnsymbol{footnote}}
\maketitle			
\renewcommand{\thefootnote}{\arabic{footnote}}

\tableofcontents

 \begin{abstract}

       A geometric  interpretation  of  curvature  and  torsion  of
   linear transports along paths is presented. A number of  (Bianchi
   type) identities satisfied by these quantities  are  derived.  The
   obtained results contain as  special  cases  the  corresponding
   classical  ones  concerning  curvature  and  torsion  of   linear
   connections.
 \end{abstract}

\section {\bfseries Introduction}
\label{I}
\setcounter {equation} {0}

     The properties of curvature and torsion tensors  of
a linear connection and the  satisfied  by  them  Bianchi
identities are well-known~\cite{Schouten/Ricci,  Helgason}.
Looking  over  the  given   in~\cite{f-LTP-Cur+Tor}
definitions and properties  of  the  curvature  and  torsion  of
linear transports along paths, one can expect to find out similar
results in this  more general case too.
To their derivation is devoted this paper.
     Sect.~\ref{II} reviews some definitions and results
from~\cite{f-LTP-general,f-LTP-Cur+Tor}
and also contains new ones needed for our investigation.
Sect.~\ref{III}  proposes a geometrical  interpretation of the torsion
of a linear transport along paths on the basis of the question of
the existence of an `infinitesimal' parallelogram. Sect.~\ref{IV} deals
with the geometrical meaning of the curvature of a linear transport
along paths.  It  is shown that the curvature governs the main change
of a vector  after a suitable transportation along a `small'
(infinitesimal) close path. Sect.~\ref{V} gives the derivation
of the generalizations of the  Bianchi
identities in the case of linear transports along paths. This is done
by using the developed in~\cite{f-Jacobi} method for obtaining
many-point generalizations of the Jacobi identity.
Sect.~\ref{VI} closes the paper with some concluding remarks,
including a criterion for flatness of a linear transport along paths.

\section {\bfseries Some preliminary definitions and results}
\label{II}
\setcounter {equation} {0}

     Below are summarized some needed for this  investigation
definitions  and
results for linear transports along paths in vector bundles and their
curvature and torsion.

     Let $(E,\pi,B)$ be a real%
\footnote{All of the results of this work are valid mutatis mutandis
in the complex case too.}
vector  bundle  with  base  $B$,  total space $E$ and
projection $\pi:E\to B$. The fibres  $\pi^{-1}(x)$,
$x\in B$ are supposed to be isomorphic real vector spaces.
Let  $\gamma:J\to B$, with $J$ being a real
interval, be an arbitrary path in $B$.
According  to~\cite[definition  2.5]{f-LTP-general} a
linear  transport (L-transport) in $(E,\pi,B)$ is a map
$L:\gamma\mapsto L^\gamma$, where the L-transport along $\gamma$ is
$L^\gamma:(s,t)\mapsto  L^\gamma_{s\to t}$, $s,t\in J$.
Here
$L^\gamma_{s\to t}: \pi^{-1}(\gamma(s))\to \pi^{-1}(\gamma(t))$
is the L-transport along $\gamma$ from $s$ to $t$.
It satisfies the equations:
\begin{eqnarray}
& & \!\!\! \!\!\! \!\!\! \!\!\! \!
L^\gamma_{s\to t}(\lambda u + \mu v) =
\lambda L^\gamma_{s\to t}u + \mu L^\gamma_{s\to t}v, \
\lambda,\mu \in {\mathbb{R}}, \  u,v \in \pi^{-1}(\gamma(s)),
\label{2.1} \\
& &
L^\gamma_{t\to r} \circ L^\gamma_{s\to t} = L^\gamma_{s\to r}, \quad
r,s,t \in J,    \label{2.2} \\
& &
L^\gamma_{s\to s} = id_{\pi^{-1}(\gamma(s))}   \label{2.3}
\end{eqnarray}
\noindent with $id_U$ being the identity map of the set~$U$.

     Propositions 2.1 and 2.3
of~\cite{f-LTP-general}  state  that  the  general
structure of $L^\gamma_{s\to t}$ is
\begin{equation}
L^\gamma_{s\to t}= \left( F^\gamma_t\right) ^{-1} \circ F^\gamma_s,
\ s,t\in J     \label{2.4}
\end{equation}
\noindent where the map $F^\gamma_s:\pi^{-1}(\gamma(s))  \to  V$  is a
linear isomorphism on a  vector  space  $V$.
The map $F^\gamma_s$ is defined up to a left composition with
a linear isomorphism $D^\gamma:V\to \underline{V}$, with
$\underline{V}$ being a vector space, i.e. up to the change
$F^\gamma_s \to D^\gamma \circ F^\gamma_s $.

     Let $\{e_i(s)\}$ be a basis  in  $\pi^{-1}(\gamma(s))$.
Here and below the indices $i,j,k,...$ run from 1 to
$\dim(\pi^{-1}(x))=\mathrm{const}=:n$.
The matrix of the L-transport $L$ (see~\cite[p.~5]{f-LTP-general}),
$H(t,s;\gamma)= \left[ H^i_j(t,s;\gamma) \right] = H^{-1}(s,t;\gamma)$,
 is defined by
$L^\gamma_{s\to t}e_j(s)=H^i_j(t,s;\gamma)e_i(t) $, where   hereafter
summation over repeated indices is assumed.
The matrix of the coefficients
of $L$~\cite[p.~13]{f-LTP-general} is
$\Gamma_\gamma(s) = \left[ \Gamma^i_j(s;\gamma) \right] =
{\partial H(s,t;\gamma)/\partial t }\left.\right|_{t=s}$.
Therefore for a $C^2$ L-transport, we have
\begin{eqnarray} & &
H^{\pm 1}(s+\varepsilon,s;\gamma) = H^{\mp 1}(s,s+\varepsilon;\gamma)=
 \nonumber \\
 & &
= \openone \mp \varepsilon \Gamma_\gamma(s) + {\varepsilon^2 \over 2}
\left(\Gamma_\gamma(s) \Gamma_\gamma(s) \mp {\partial \Gamma_\gamma(s)
\over \partial s} \right) + O(\varepsilon ^2),  \label{2.5}
\end{eqnarray}
\noindent
 with $\openone$ being the unit matrix. Here we have used
		\begin{equation} \label{2.6}
		\begin{array}{c}
\rule[-1.123456789em]{0em}{0em}
\left. {\partial^2 H(s,t;\gamma) \over \partial t^2}
\right|_{t=s} =
\Gamma_\gamma(s) \Gamma_\gamma(s) +
{\partial \Gamma_\gamma(s) \over \partial s},
\\
\left. {\partial^2 H(t,s;\gamma) \over \partial t^2}\right|_{t=s} =
\Gamma_\gamma(s) \Gamma_\gamma(s) -
{\partial \Gamma_\gamma(s) \over \partial s}.
		\end{array}
		\end{equation}
\noindent
These equations follow from the fact that the general form of the
matrix $H$ is
$H(t,s;\gamma)=F^{-1}(t;\gamma) F(s;\gamma)$ for some nondegenerate
matrix  function $F$~\cite{f-LTP-general}.

     Let $\eta:J\times J^\prime\to M$, $J$, with $J^\prime$ being
${\mathbb R} $-intervals, be a $C^2$ map on the real
differentiable manifold $M$
with a tangent bundle $(T(M),\pi,M)$.
Let $\eta(\cdot,t):s\mapsto \eta(s,t)$
and $\eta(s,\cdot):t\mapsto \eta(s,t)$, $(s,t)\in J\times J^\prime$.
Here by
$\eta^\prime(\cdot,t)$ and $\eta^{\prime\prime}(s,\cdot)$
we denote the tangent to
$\eta(\cdot,t)$ and $\eta(s,\cdot)$, respectively, vector fields.

     By~\cite[definition 2.1]{f-LTP-Cur+Tor} the torsion (operator)
of  a $C^1$ L-transport $L$ in $(T(M),\pi,M)$ is a map
\[
{\mathcal T}:\eta\mapsto{\mathcal T}^{\,\eta}:J\times J^\prime\to T(M)
\]
such that
\begin{equation}
{\mathcal T}^{\,\eta} (s,t):= {\mathcal D}^{\eta(\cdot,t)}_s
\eta^{\prime\prime}(\cdot,t) -
{\mathcal D}^{\eta(s,\cdot)}_t \eta^{\prime}(s,\cdot)
					\in T_{\eta(s,t)}(M),
							\label{2.7}
\end{equation}
\noindent
where ${\mathcal D}^\gamma_s$ is the associated with
$L$ differentiation along paths~\cite{f-LTP-general}, defined by
\[
{\mathcal D}^\gamma_s \sigma :=
\left( {\mathcal D}^\gamma \sigma \right) (\gamma(s)) :=
\left. \left[ {\partial\over \partial \varepsilon}\left(
L^\gamma_{s+\varepsilon \to s}\sigma(s+\varepsilon) \right) \right]
\right|_{\varepsilon=0}
\]
for a $C^1$ section $\sigma$.

     Analogously~\cite{f-LTP-Cur+Tor}, for
$\eta:J\times  J^\prime \to B$ the
curvature (operator) of an \mbox{L-tran}sport $L$ in the vector bundle
$(E,\pi,B)$ is a map
\[
{\mathcal R}: \eta\mapsto {\mathcal R}^\eta:(s,t) \mapsto
{\mathcal R}^\eta(s,t):\mathrm{Sec}^2(E,\pi,B)\to\mathrm{Sec}(E,\pi,B)
\]
such that
\begin{equation}   \label{2.8}
{\mathcal R}^\eta(s,t):=
{\mathcal D}^{\eta(\cdot,t)}\circ {\mathcal D}^{\eta(s,\cdot)} -
{\mathcal D}^{\eta(s,\cdot)}\circ {\mathcal D}^{\eta(\cdot,t)}.
\end{equation}

	In terms of the coefficient matrix $\Gamma$ the components of
torsion and curvature are respectively~\cite{f-LTP-Cur+Tor}
\begin{eqnarray}
& &
\left({\mathcal T}^{\,\eta}(s,t) \right)^i =
\Gamma^i_j(s;\eta(\cdot,t)) \left( \eta^{\prime\prime}(s,t)\right)^j -
\Gamma^i_j(t;\eta(s,\cdot)) \left( \eta^{\prime}(s,t) \right)^j ,
						\label{2.9} \\
& &
\left[ \left( {\mathcal R}^\eta(s,t)\right) ^i_j \right]
= {\partial\over\partial s}\Gamma_{\eta(s,\cdot)}(t) -
{\partial\over{\partial t}} \Gamma_{\eta(\cdot,t)}(s) \> + \nonumber\\
 & & \hspace{24.3mm}
+ \> \Gamma_{\eta(\cdot,t)}(s)\Gamma_{\eta(s,\cdot)}(t) -
\Gamma_{\eta(s,\cdot)}(t)\Gamma_{\eta(\cdot,t)}(s). 	\label{2.10}
\end{eqnarray}

	Below we shall need the following definitions:
\begin{defi} \label{d2.1}
The torsion vector field (operator) of an L-transport in the tangent to
a manifold bundle  is a section $T^{\,\eta}\in$\
\(
\mathrm{Sec}
\left(\left.\left(T(M),\pi,M\right)\right|_{\eta(J,J')}\right)
\)
defined by
\begin{equation}
 T^{\,\eta}(\eta(s,t)):={\mathcal T}^{\,\eta}(s,t).	\label{2.11}
\end{equation}
\end{defi}
Defining
\(
({\mathcal D}^\gamma\sigma)(\gamma(s)) :=
			{\mathcal D}_{s}^{\gamma}\sigma,
\)
from~(\ref{2.7}) we get
	\begin{equation}  \label{2.12}
T^{\,\eta}(\eta(s,t)) = \left(
{\mathcal D}^{\eta(\cdot,t)}\eta ''(\cdot,t) -
{\mathcal D}^{\eta(s,\cdot)}\eta '(s,\cdot)
\right)(\eta(s,t)).
	\end{equation}

\begin{defi} \label{d2.2}
The curvature vector field (operator) of an L-transport
is a $C^2$ section
\(
R^\eta\in
\mathrm{Sec}^2\left( \left. (E,\pi,B) \right|_{\eta(J,J')}\right)
\)
defined by
\begin{equation}
R^\eta(\eta(s,t)):={\mathcal R}^\eta(s,t).		\label{2.13}
\end{equation}
\end{defi}

\begin{defi} \label{d2.3}
An L-transport along paths is called flat ($\equiv$curvature free)
on a set $U\subseteq B$ if its curvature operator vanishes on~$U$.
It is called flat if it is flat on $B$, i.e. in the case $U=B$.
\end{defi}

\section {\bfseries Geometrical interpretation of the torsion}
\label{III}
\setcounter {equation} {0}

	Let $\eta:J\times J'\to M$ be a $C^1$ map into the manifold $M$,
$(s,t)\in J\times J'$, and $\delta,\varepsilon\in{\mathbb{R}}$ be
such that $(s+\delta,t+\varepsilon)\in J\times J'$. Below we consider
$\delta$ and $\varepsilon$ as `small' (infinitesimal) parameters with
respect to which expansions like (\ref{2.5}) will be used.

\begin{figure}							  
\vspace*{-1.8em}						  
{\tiny      \input{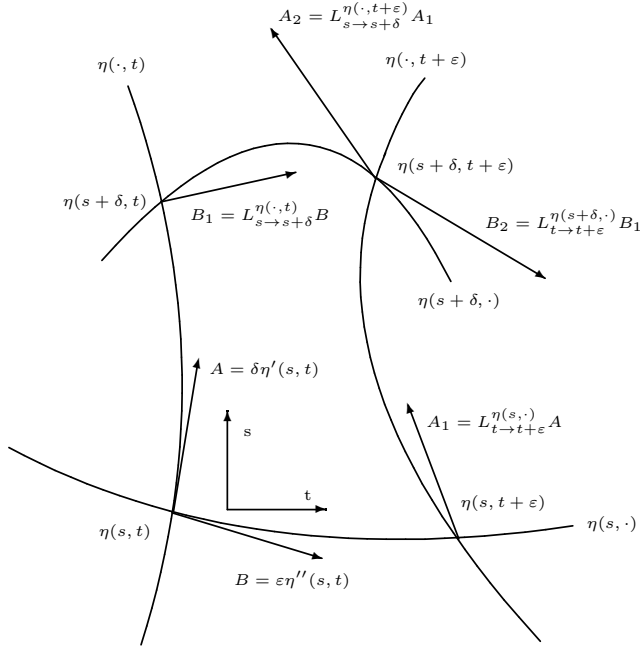}	}			  
\caption{Geometrical interpretation of the torsion} \label{Fig1}  
\end{figure}

	Consider the following two paths from $\eta(s,t)$ to
$\eta(s+\delta,t+\varepsilon)$ (see figure~\ref{Fig1}):
the first, through $\eta(s+\delta,t)$,
being a product of
$\eta(\cdot,t):[s,s+\delta]\to M$
and
$\eta(s+\delta,\cdot):[t,t+\varepsilon]\to M$,
and the second one, through
$\eta(s,t+\varepsilon)$,
being a product of
$\eta(s,\cdot):[t,t+\varepsilon]\to M$
and
$\eta(\cdot,t+\varepsilon):[s,s+\delta]\to M$.
(Here $\delta$ and $\varepsilon$ are considered as positive, but this
is inessential.)


	Up to $O(\delta^2)$  and   $O(\varepsilon^2)$     the
vectors $A:=\delta\eta'(s,t)$ and  $B:=\varepsilon\eta''(s,t)$ are the
displacement  vectors~\cite{f-DE-TP}
(linear elements~\cite{Schouten/Ricci}),
respectively,  of  $\eta(s+\delta,t)$  and  $ \eta(s,t+\varepsilon)$
with respect to $\eta(s,t)$.

	Using (\ref{2.5}) and keeping only the first order in
$\varepsilon$ and $\delta$ terms in it, we get the following
\textsl{component} relation:
\begin{equation}
\left( L_{s\to s+\delta}^{\eta(\cdot,t)} B \right)^i -
\left( L_{t\to t+\varepsilon}^{\eta(s,\cdot)} A \right)^i =
(B-A)^i -
\delta\varepsilon\left( {\mathcal T}^{\,\eta}(s,t) \right)^i +
O(\delta\varepsilon^2) + O(\delta^2\varepsilon).  \label{3.1}
\end{equation}

	According to~\cite[ch. V, sect. 1]{Schouten/physics} this
result has the following interpretation.
\textit{%
After the `L-transportation' of two  linear elements  $A$ and
$B$ along each other we get, up to second order terms,
a pentagon with a closure vector
$-\delta\varepsilon{\mathcal T}^{\,\eta}(s,t)$.
	}%
This implies the existence of an infinitesimal parallelogram
only in the torsion free case.

	Using again~(\ref{2.5}) and keeping only first order terms,
after some algebra, we find
\begin{eqnarray}
& &  \nonumber
\left(
L_{t\to t+\varepsilon}^{\eta(s+\delta,\cdot)} \circ
			 L_{s\to s+\delta}^{\eta(\cdot,t)} B -
L_{s\to s+\delta}^{\eta(\cdot,t+\varepsilon)} \circ
			L_{t\to t+\varepsilon}^{\eta(s,\cdot)} A
\right)^i
 \> = \\
& & \nonumber
 = \>
\left[
\left( L_{t\to t+\varepsilon}^{\eta(s,\cdot)} B \right)^i -
\left( L_{s\to s+\delta}^{\eta(\cdot,t)} A \right)^i
\right] -
    \delta\varepsilon\left({\mathcal T}^{\,\eta}(s,t)\right)^i
 \> + \\
 & & 	 \label{3.2}
 \quad\> + \>
   O(\delta^3) + O(\delta^2\varepsilon) + O(\delta\varepsilon^2) +
   O(\varepsilon^3).
\end{eqnarray}

	Note that if $\eta$ is a family of L-paths, i.e.
$ L_{s_1\to s_2}^{\eta(\cdot,t)}\eta'(s_1,t)=\eta^\prime(s_2,t) $ and
\(
L_{t_1\to t_2}^{\eta(s,\cdot)}\eta''(s,t_1) =
				\eta{^\prime}{^\prime}(s,t_2),
\)
for all $s,s_1,s_2\in J$ and $t,t_1,t_2\in J^\prime$, the expression
in the square brackets in~(\ref{3.2}) is simply $(B-A)^i$.

	So, the torsion describes the first order correction to the
difference of two (infinitesimal) displacement vectors when they are
(L-)transported in the above-described way.

\section {\bfseries Geometrical interpretation of curvature}

\label{IV}
\setcounter {equation} {0}

	Let $(E,\pi,B)$ be a vector bundle,
$\eta:J\times J^\prime\to B$ be a
$C^1$ map, and $L$ be a $C^2$ L-transport along paths in $(E,\pi,B)$.
Let
$(s,t)\in J\times J^\prime$ and
$\delta, \varepsilon \in {\mathbb{R}}$, $\delta,\varepsilon >0$
(this condition is insignificant for the final result)
be such that $(s+\delta,t+\varepsilon)\in J\times J^\prime$.

\begin{figure}							   
{\scriptsize	\input{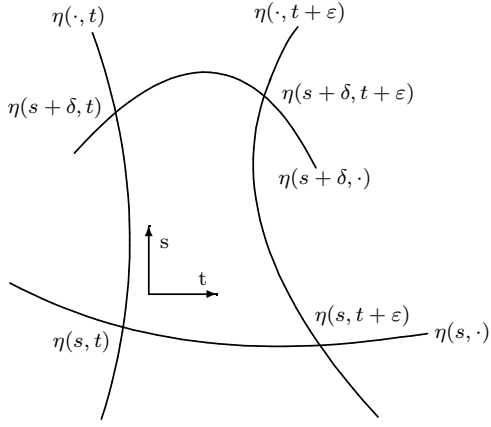}	}			   
\caption{Geometrical interpretation of the curvature} \label{Fig2} 
\end{figure}							   

	Consider the paths on figure~\ref{Fig2}.
	The result of an L-transport of a vector
from $\eta(s,t)$ to $\eta(s+\delta,t)$ along
$\left.\eta(\cdot,t)\right|_{[s,s+\delta]}$,
then from $\eta(s+\delta,t)$ to $\eta(s+\delta,t+\varepsilon)$  along
$\left.\eta(s+\delta,\cdot)\right|_{[t,t+\varepsilon]}$,
then from $\eta(s+\delta,t+\varepsilon)$ to $\eta(s,t+\varepsilon)$
along $\left.\eta(\cdot,t+\varepsilon)\right|_{[s,s+\delta]}$,
and, at last, from $\eta(s,t+\varepsilon)$ to $\eta(s,t)$ along
$\left.\eta(s,\cdot)\right|_{[t,t+\varepsilon]}$
is expressed by
		\begin{prop}		\label{p4.1}
	For any $C^2$ L-transport, we have
		\begin{eqnarray}
& &   \!\!\!\! \!\!\!\! \!\!\!\!\!
L_{t+\varepsilon\to t}^{\eta(s,\cdot)} \circ
L_{s+\delta\to s}^{\eta(\cdot,t+\varepsilon)} \circ
L_{t\to t+\varepsilon}^{\eta(s+\delta,\cdot)} \circ
L_{s\to s+\delta}^{\eta(\cdot,t)} =      \nonumber \\
& &
\!\!\!\!\!\!\! \!\!\!\!
=id_{\pi^{-1}(\eta(s,t))} -
\delta\varepsilon{\mathcal R}^{\eta}(s,t) +
O(\delta^3) + O(\delta^2\varepsilon) +O(\delta\varepsilon^2) +
O(\varepsilon^3). \label{4.1}
		\end{eqnarray}
		\end{prop}

	\textit{Proof.}
In a field $ \{e_i(s,t),\ (s,t)\in J\times J^\prime\}$
of bases in $\pi^{-1}(\eta(J,J^\prime))$ the matrix of the linear map
standing in the l.h.s. of (\ref{4.1}) is
$$
H(t,t+\varepsilon;\eta(s,\cdot))
H(s,s+\delta;\eta(\cdot,t+\varepsilon))
H(t+\varepsilon,t;\eta(s+\delta,\cdot))
H(s+\delta,s;\eta(\cdot,t)).
$$
Substituting here (\ref{2.5}) and using the expressions
\[
\Gamma_{\eta(s+\delta,\cdot)}(t)=\Gamma_{\eta(s,\cdot)}(t) +
\delta{\partial\over\partial s}\Gamma_{\eta(s,\cdot)}(t) +O(\delta^2),
\]
$$
\Gamma_{\eta(\cdot,t+\varepsilon)}(s)=\Gamma_{\eta(\cdot,t)}(s) +
\varepsilon{\partial\over\partial t}\Gamma_{\eta(\cdot,t)}(s) +
O(\varepsilon^2),
$$
we find this matrix to be
	\begin{eqnarray*}
&\!\!\!&
\!\!\!\!\!\!\!\!\!\!\!\openone \! + \!  \delta\varepsilon\left(
{\partial\over{\partial t}}\Gamma_{\eta(\cdot,t)}(s) -
{\partial\over{\partial s}}\Gamma_{\eta(s,\cdot)}(t) -
\Gamma_{\eta(\cdot,t)}(s) \Gamma_{\eta(s,\cdot)}(t) +
\Gamma_{\eta(s,\cdot)}(t) \Gamma_{\eta(\cdot,t)}(s) \right) \!\!+ \\
& &
+ \> O(\delta^3) +O(\delta^2\varepsilon) +O(\delta\varepsilon^2)
+O(\varepsilon^3).
	\end{eqnarray*}
Taking into account (\ref{2.10}),
we get this expression as
	\begin{equation} \label{4.1a}
\openone - \delta\varepsilon
\left[ \left({\mathcal R}^\eta(s,t)\right)_{j}^{i} \right]
+O(\delta^3)+O(\delta^2\varepsilon) +O(\delta\varepsilon^2) +
						O(\varepsilon^3)
	\end{equation}
which is
simply the matrix form of (\ref{4.1}).~\qed

	Proposition~\ref{p4.1} shows that
\textit{%
up to third order terms the result of the above-described
transportation of a vector $A\in\pi^{-1}(\eta(s,t))$ is%
}
\begin{equation}	\label{4.2}
A -\delta\varepsilon\left( {\mathcal R}^\eta(s,t)\right) (A).
\end{equation}

	Another corollary of (\ref{4.1}) is the (equivalent to the
definition of the curvature) equality
	\begin{equation}	 \label{4.3}
{\mathcal R}^\eta(s,t) = -
\lim_{ \stackrel{\delta\to 0}{\varepsilon\to 0} }
\left[{1\over \delta\varepsilon}
\left(
L_{t+\varepsilon\to t}^{\eta(s,\cdot)}\!\circ\!
L_{s+\delta\to s}^{\eta(\cdot,t+\varepsilon)}\!\circ\!
L_{t\to t+\varepsilon}^{\eta(s+\delta,\cdot)}\!\circ\!
L_{s\to s+\delta}^{\eta(\cdot,t)} -
id_{\pi^{-1}(\eta(s,t))}
\right)
\right].
	\end{equation}

\section {\bfseries Bianchi-type identities}

\label{V}
\setcounter {equation} {0}

       	The curvature operator (\ref{2.8}) is simply a commutator of
two derivations along paths. As we shall see below,
the torsion (\ref{2.7}) is also a skewsymmetric expression.
All this allows one to apply the developed in~\cite{f-Jacobi} method
for obtaining Jacobi-type identities. This can be done as follows.

	Let us take an arbitrary map $\tau^k:J^k\to B$, with
$J^k:=J\times\cdots\times J$, where $J$ appears $k$ times,
$k\in \mathbb{N}$,
and $B$ being the base of the vector bundle $(E,\pi,B)$. Let
$s:=(s^1,...,s^k)\in J^k$. We define the $C^1$ paths
$\tau_a:J\to B$ by
$\tau_a(\sigma):=\left.\tau^k(s)\right|_{s^a=\sigma}$,
$\sigma\in J$ and the
maps (families of paths) $\tau_{ab}:J\times J\to B$ by
$\tau_{ab}(\sigma_1,\sigma_2):=\left. \tau^k(s)\right|_{s^a=\sigma_1,
s^b=\sigma_2}$, $ \sigma_1,\sigma_2\in J$
which depend implicitly on s.
Hereafter
$a,b,c,d=1,\ldots,k$.  We write  $\dot\tau_a$ for the tangent to
$\tau_a$ vector field in the case when $ (E,\pi,B)=(T(M),\pi,M) $
for some manifold $M$.

	\begin{prop}\label{p5.1}
The following properties of antisymmetry are valid:
	\begin{eqnarray}
& &
{\mathcal R}^{\tau_{ab}}(s_a,s_b) +
		{\mathcal R}^{\tau_{ba}}(s_b,s_a)=0,
						\label{5.1}\\ & &
T^{\tau_{ab}} + T^{\tau_{ba}} =0
\ {\mathrm{or}}\ {\mathcal T}^{\,\tau_{ab}}(s_a,s_b) +
{\mathcal T}^{\,\tau_{ba}}(s_b,s_a) =0.		\label{5.2}
	\end{eqnarray}
	\end{prop}

	{\bfseries Remark}
These equalities are analogues of the usual skewsymmetry of curvature
and torsion tensors in the tensor analysis~\cite{Schouten/Ricci}.

   \textit{Proof.} The two-point Jacobi-type identity is
$ \left( (A_{ab})_{[a,b]} \right)_{<a,b>} \equiv 0  $
(see~\cite[eq. (5.1)]{f-Jacobi}), where  $ A_{ab}  $
are elements of an Abelian group,
$ (A_{ab})_{[a,b]}:=A_{ab} - A_{ba} $
and
$ \left(A_{ab}\right)_{<a,b>}:=A_{ab} + A_{ba}  $.
Substituting here
$ A_{ab}={\mathcal D}^{\tau_a} \circ {\mathcal D}^{\tau_b} $
in the case of a vector bundle  $ (E,\pi,B)  $ and
$ A_{ab}={\mathcal D}^{\tau_a}({\dot\tau_b})$
in the case of the tangent to a manifold $M$ bundle  $(T(M),\pi,M)$
and using (\ref{2.8}) and (\ref{2.12}) (or (\ref{2.7})), one gets
respectively (\ref{5.1}) and (\ref{5.2}).~\qed

\begin{prop} \label{p5.2}
	The following identities are valid:
	\begin{eqnarray}
& &
	\begin{array}{l}
\left\{
{\mathcal D}^{\tau_a} \circ {\mathcal R}^{\tau_{bc}}(s_b,s_c) -
{\mathcal R}^{\tau_{bc}}(s_b,s_c) \circ {\mathcal D}^{\tau_a}
\right\} _{<a,b,c>} \equiv 0 \ \mathrm{or}\
 \\
\left\{
{\mathcal D}^{\tau_a} \left( {\mathcal R}^{\tau_{bc}}(s_b,s_c)
\right)
\right\}_{<a,b,c>} \equiv 0,
	\end{array}
\label{5.3} \\
& & \
\left\{
\left( {\mathcal R}^{\tau_{ab}}(s_a,s_b) \right)(\dot\tau_c)
\right\}
_{<a,b,c>} \equiv \left\{ {\mathcal D}^{\tau_a}
\left( {\mathcal T}^{\,\tau_{bc}} \right) \right\}_{<a,b,c>},
							\label{5.4}
\end{eqnarray}
\end{prop}
\noindent
where $<\ldots>$ means summation over the cyclic permutations of the
corresponding indices.

	{\bfseries Remark.}
These identities are analogous, respectively, of
the second and first Bianchi identities in tensor
analysis~\cite{Schouten/physics,Helgason}.
This is clear from the fact that due to the antisymmetries
(\ref{5.1}) and (\ref{5.2}) the cyclization over the
indices $a,\ b\  \mathrm{and}\ c$, i.e. the operation $<\ldots>$,
in (\ref{5.3}) and (\ref{5.4}) may be replaced with
antisymmetrization over the indices $a,\ b \ \mathrm{and}\ c$.
(E.g. if
$A_{abc}=-A_{acb}$ and
\(
\left(A_{abc}\right)_{[a,b,c]} :=
\left(A_{abc}+A_{bca}+A_{cab}\right)_{[b,c]},
\)
then
$2\left(A_{abc}\right)_{<abc>} = \left(A_{abc}\right)_{[abc]}$.)

	\textit{Proof.} The (3-point) generalized Jacobi identity
(see~\cite[eq. (5.2)]{f-Jacobi}) is
$\left(\left(A_{abc}\right)_{[a,[b,c]]}\right)_{<a,b,c>} \equiv 0,$
with $A_{abc}$ being elements of an Abelian group,
\(
\left(A_{abc}\right)_{[a,[b,c]]} :=
\left(A_{abc} - A_{bca}\right)_{[b,c]}
\)
and
$\left(A_{abc}\right)_{<a,b,c>}:=A_{abc} + A_{bca} + A_{cab}$.

	We put
$A_{abc}=\mathcal{D}^{\tau_a}\circ\mathcal{D}^{\tau_b}\circ
\mathcal{D}^{\tau_c} $
in the vector bundle case and
$A_{abc}=\left(\mathcal{D}^{\tau_a} \circ
\mathcal{D}^{\tau_b}\right) (\dot\tau_c)$
in the tangent bundle case. In this way, after some simple algebra
(see~(\ref{2.8}), (\ref{2.7}) and (\ref{2.1})-(\ref{2.3})), we get
respectively (\ref{5.3}) and (\ref{5.4}).~\qed

	The 4-point generalized Jacobi-type identity
\[
\left\{ \left( A_{abcd} \right) _ {[a,[b,[c,d]]]} +
\left( A_{adcb} \right) _ {[a,[d,[c,b]]]} \right\}_{<a,b,c,d>}\equiv 0
\]
with
$ \left( A_{abcd} \right) _ {[a,[b,[c,d]]]} :=
\left( A_{abcd} - A_{bcda} \right) _ {[b,[c,d]]} $
and
$ \left( A_{abcd} \right) _ {<a,b,c,d>} :=
 A_{abcd} + A_{bcda} + A_{cdab} + A_{dabc}$
also produces an interesting identity in our case. In fact, putting
$  A_{abcd} =
\mathcal{D}^{\tau_a}\circ\mathcal{D}^{\tau_b}\circ
\mathcal{D}^{\tau_c}\circ\mathcal{D}^{\tau_d}$
in the vector bundle case, one can easily prove after some simple
calculations

\begin{prop} \label{p5.3}
The identity
\begin{equation}
\left\{ \mathcal{R}^{\tau_{ab}}(s_a,s_b) \left(R^{\tau_{cd}} \right)
\right\}_{<a,b,c,d>} \equiv 0,			\label{5.5}
\end{equation}
where $R^{\tau_{cd}}$ is the curvature vector field on
$\tau^k(J,\ldots,J)$ is valid.
\end{prop}

	{\bfseries Remark.} This result generalizes eq. (6.5)
of~\cite{f-Jacobi} in the classical tensor case.

	The last result also follows from the evident chain identity
	\begin{eqnarray*}
0 &\equiv&
\left\{
\mathcal{R}^{\tau_{ab}}(s_a,s_b)
			\circ \mathcal{R}^{\tau_{cd}}(s_c,s_d) -
\mathcal{R}^{\tau_{ab}}(s_a,s_b)
			\circ \mathcal{R}^{\tau_{cd}}(s_c,s_d)
\right\}_{<a,b,c,d>} \equiv   \\
  &\equiv&
\left\{
\mathcal{R}^{\tau_{ab}}(s_a,s_b)
			\circ \mathcal{R}^{\tau_{cd}}(s_c,s_d) -
\mathcal{R}^{\tau_{cd}}(s_c,s_d)
			\circ \mathcal{R}^{\tau_{ab}}(s_a,s_b)
\right\}_{<a,b,c,d>} \equiv    \\
  &\equiv&
\left\{
\left(
\mathcal{R}^{\tau_{ab}}(s_a,s_b) \left( R^{\tau_{cd}} \right)
\right)
\left( \tau_{cd}(s_c,s_d) \right)
\right\}_{<a,b,c,d>} \equiv    \\
  &\equiv&
\left( \left\{
\mathcal{R}^{\tau_{ab}}(s_a,s_b) \left( R^{\tau_{cd}} \right)
\right\}_{<a,b,c,d>} \right) \left(\tau^k(s)\right).
	\end{eqnarray*}

	Note that in the tangent bundle case the substitution
\[
A_{abcd} = \left(
\mathcal{D}^{\tau_a}\circ\mathcal{D}^{\tau_b}\circ
\mathcal{D}^{\tau_c} \right) \bigl( \dot\tau_d \bigr)
\]
leads to the trivial identity $0\equiv 0$.

\section {\bfseries Conclusion}
\label{VI}
\setcounter {equation} {0}

	In this paper we have examined some natural properties of the
curvature (resp. the torsion) of linear transports along paths in
vector bundles (resp. in the tangent bundle to a manifold).
These properties are similar to the ones in the theory of linear
connections. The cause for this similarity is that in the case of the
parallel transport assigned to a linear connection our results
reproduce the corresponding ones in the classical tensor analysis.
The reduction to the known classical results can easily be proved
by applying the used in~\cite{f-LTP-Cur+Tor}  method for introduction
of curvature and torsion of a linear connection by means of its
parallel transport.

	In connection with this below is presented the generalization
of the theorem that a linear connection is flat iff the assigned to it
parallel transport is independent of the path (curve) along which it
acts and depends only on the initial and final points of
the transportation.

	\begin{thm} \label{th6.1}
An L-transport in  $ (E,\pi,B) $  is flat on  $ U\subseteq B $
if and only if in~$U$ it is independent of the path
(lying in~$U$) along which it acts and depends only on its initial
and final points, i.e. the set $\{ L_{s\to t}^{\gamma}\}$ forms a flat
L-transport in  $U\subseteq B$ iff $L_{s\to t}^{\gamma}$ for
$\gamma:J\to U$ depends only on  the points $\gamma(s)$ and
$\gamma(t)$, but not on the path $\gamma$ itself.
	\end{thm}

	{\bfseries Remark.}
In this theorem we implicitly suppose~$U$ to be linearly connected,
i.e. its every two points can be connected by a path
lying entirely in~$U$. Otherwise the theorem may not be true.

	\textit{Proof.} Let the L-transport  $L$ be flat, i.e.
\(
\mathcal{R}^\eta(s,t)\equiv 0 \ \mathrm{for}\
\eta:J\times J^\prime\to U\subseteq B.
\)
By~\cite[theorem 3.1]{f-LTP-Cur+Tor} there is a field of bases
$\{e_i\}$ on~$U$ in which the matrix of  $L$ is unit, i.e.
$H(t,s;\gamma)=\openone$,  $\gamma:J\to U$. In these bases for
$u\in\pi^{-1}(\gamma(s))$, we have
$L_{s\to t}^{\gamma}u=
H_j^i(t,s;\gamma)u^j\left(e_i\left.\right|_{\gamma(t)}\right)=
u^i\left(e_i\left.\right|_{\gamma(t)}\right)$,
which evidently depends on the points  $\gamma(s)$ and $\gamma(t)$
but not on the path  $\gamma$ itself.

	Conversely, let for  $\gamma:J\to U$  the transport
$L_{s\to t}^{\gamma}$ depends only on the points
$\gamma(s)$ and $\gamma(t)$ and not on the path  $\gamma$
connecting them. For fixed  $x_0\in U$ and basis
$\{ e_{i}^{0} \}$ in $\pi^{-1}(x)$
we define on~$U$ the field of bases  $\{e_i\}$ by
$e_i\left.\right|_x := L_{a\to b}^{\beta} e_{i}^{0}$,
where  $\beta$ is
any path in~$U$ joining  $x_0$ and  $x\in U$, and such that
$\beta(a):=x_0\ \mathrm{and}\ \beta(b):=x$.
By assumption $\{ e_i\left.\right|_x \}$ depends
only on  $x$ but not on  $\beta$. Using that $L_{s\to t}^{\gamma}$
depends only on $\gamma(s)$ and $\gamma(t)$, we have
	\begin{eqnarray*}
L_{s\to t}^{\gamma} \left( { e_i\left.\right|_{\gamma(s)} }\right) &=&
L_{a\to b}^{\alpha} \left( { e_i\left.\right|_{\alpha(a)} }\right)\>=\\
&=&\>
L_{a\to b}^{\alpha} \left( L_{c\to a}^{\alpha} { e_i^0} \right)   =
L_{c\to b}^{\alpha} { e_i^0} = { e_i\left.\right|_{\alpha(b)} }   =
{ e_i\left.\right|_{\gamma(t)} },
	\end{eqnarray*}
\noindent where  $\alpha$ is any path in~$U$  such that
$\alpha(a)=\gamma(s),\ \alpha(b)=\gamma(t)$, and  $\alpha(c)=x_0$.
As
$L_{s\to t}^{\gamma} \left( { e_i\left.\right|_{\gamma(s)} } \right)=
H_{i}^{j}(t,s;\gamma) { e_j\left.\right|_{\gamma(t)} }$,
we see that in
$\{e_i\}$ the matrix of  $L$ is  $H(t,s;\gamma)=\openone$,
which, again by~\cite[theorem 3.1]{f-LTP-Cur+Tor},
implies the flatness of $L$ in~$U$.~\qed

	In conclusion we have to note that all of the results of the
present paper remain true in the complex case. For this purpose one has
simply to replace in it the word `real' with `complex' and the symbols
$\mathbb{R}$ and $\dim$ with $\mathbb{C}$  and  $\dim_\mathbb{C}$
respectively.

\section * {\bfseries Acknowledgement}
\label{VII}
\setcounter {equation} {0}

This work was partially supported by the National Science Foundation
of Bulgaria under Grant No. F642.

\bibliography{bozhopub,bozhoref}
\bibliographystyle{unsrt}

\end{document}